\documentclass[12pt]{amsart}
\usepackage[]{amsmath, amsthm, amsfonts, verbatim, amssymb}
\usepackage{enumerate}
\newcommand\C{{\mathbb C}}

\newcommand\dee{\partial}
\newcommand\Om{\Omega}
\newcommand\Obar{\overline{\Omega}}

\newcommand\Ohat{\widehat\Omega}

\numberwithin{equation}{section}

\begin{document}

\title[Szeg\H o coordinates]
{Szeg\H o coordinates, quadrature domains, and double quadrature domains}
\author[Bell, Gustafsson, Sylvan]
{Steven R.~Bell, Bj\"orn Gustafsson, and Zachary A.~Sylvan}

\address[]{Mathematics Department, Purdue University, West Lafayette,
IN  47907}
\email{bell@math.purdue.edu}

\address[]{Department of Mathematics, KTH, 100 44 Stockholm, Sweden}
\email{gbjorn@kth.se}

\address[]{Mathematics Department,
University of California
Berkeley, CA 94720}
\email{zack@math.berkeley.edu}



\begin{abstract}
We define Szeg\H o coordinates on a finitely
connected smoothly bounded planar domain which effect
a holomorphic change of coordinates on the domain
that can be as close to the identity as desired and
which convert the domain to a quadrature domain with
respect to boundary arc length.  When these Szeg\H o
coordinates coincide with Bergman coordinates, the
result is a double quadrature domain with respect to
both area and arc length.  We enumerate a host of
interesting and useful properties that such double
quadrature domains possess, and we show that such
domains are in fact dense in the realm of bounded 
$C^\infty$-smooth finitely connected domains.
\end{abstract}

\maketitle

\theoremstyle{plain}

\newtheorem {thm}{Theorem}[section]
\newtheorem {lem}[thm]{Lemma}

\hyphenation{bi-hol-o-mor-phic}
\hyphenation{hol-o-mor-phic}

\section{Introduction}
\label{intro}
The unit disc $D_1(0)$ is the quintessential example of a double quadrature
domain.  Holomorphic functions $h$ on the disc satisfy two quadrature identities:
$$\pi h(0)=\iint_{D_1(0)} h\ dA\qquad\text{and}\qquad 2\pi h(0)=\int_{bD_1(0)} h(z)\ ds,$$
whenever these integrals make sense.  We will show in this paper that double quadrature
domains are rather commonplace among smoothly bounded domains in the plane.

A bounded domain $\Om$ in the complex plane is a quadrature
domain with respect to area measure if there exist finitely many points
$\{w_j\}_{j=1}^N$ in the domain and non-negative integers $n_j$ such
that complex numbers $c_{jk}$ exist satisfying
$$\iint_\Om h\ dA = \sum_{j=1}^N\sum_{k=0}^{n_j} c_{jk} h^{(k)}(w_j)$$
for every function $h$ in the Bergman space $H^2(\Om)$ of square integrable
holomorphic functions on $\Om$.  Here, $dA$ denotes Lebesgue area
measure.  We shall say that such a domain is an {\it area quadrature
domain}.  Aharonov and Shapiro \cite{AS} proved that such domains
have piecewise real-analytic boundaries (see also \cite{G}).

A bounded domain $\Om$ bounded by finitely many non-intersecting piecewise
$C^1$-smooth Jordan curves is a quadrature domain with respect to boundary
arc length measure if there exist finitely many points $\{w_j\}_{j=1}^N$ in
the domain and non-negative integers $n_j$ such that complex numbers $c_{jk}$
exist satisfying
$$\int_{b\Om} h\ ds = \sum_{j=1}^N\sum_{k=0}^{n_j} c_{jk} h^{(k)}(w_j)$$
for every function $h$ in the Hardy space $H^2(b\Om)$ of holomorphic
functions with $L^2$ boundary values.  Here, $ds$ denotes the boundary arc
length measure.  We shall say that such a domain is a {\it boundary arc
length quadrature domain}.  Boundary arc length quadrature domains
have $C^\infty$-smooth real-analytic boundaries (see \cite{G2}, where it is
also shown that the boundary regularity assumptions can be greatly relaxed).

After some pioneering work of Avci \cite{Av} and Shapiro and Ullemar
\cite{SU}, the second author \cite{G2} studied planar quadrature domains
with respect to boundary arc length from the point of view of
half-order differentials and Riemann surface theory; he characterized
such domains and showed that they are dense in the category of bounded
domains in the plane bounded by finitely many non-intersecting
Jordan curves.  Here, we shall reframe these results in more
elementary terms and explain how to view boundary arc length quadrature
domains in terms of Szeg\H o coordinates, which are the
analogue of Bergman coordinates as developed in \cite{B4,B6,B7}.
Because our results are most interesting in the category of
$C^\infty$-smoothly bounded domains, we shall take the royal
road and restrict our attention to domains of this form.

For more information about the state of research on quadrature
domains, see \cite{GS} and the volume \cite{EGKP} that holds it,
and \cite{S}.  We also remark that it is possible to study
quadrature domains with respect to harmonic functions in the
plane, and to prove density results in this context (see
\cite{G3} and Sakai \cite{Sa}).

Suppose that $\Om$ is a bounded domain bounded by
$n$ non-intersecting $C^\infty$-smooth Jordan curves.  Let
$A^\infty(\Om)$ denote the subspace of $C^\infty(\Obar)$ consisting
of holomorphic functions.  Let
$K(z,w)$ denote the Bergman kernel associated to $\Om$ and
let $S(z,w)$ denote the Szeg\H o kernel.  Let $K^0(z,w)$ also
denote $K(z,w)$ and let $K^m(z,w)=\frac{\dee^m}{\dee\bar w^m}K(z,w)$.
As in \cite{B7}, we define the {\it Bergman span\/} associated to
$\Om$ to be the complex linear span $\mathcal K$ of all functions
$h(z)$ of $z$ of the form $K^m(z,a)$ as $a$ ranges over $\Om$ and
$m$ ranges over all non-negative integers.  If $U$ is an open
subset of $\Om$, let ${\mathcal K}_U$ denote the complex linear
span of functions of $z$ from $\{K(z,a):a\in U\}$,
and given a point $a$ in $\Om$, let ${\mathcal K}_a$ denote
the complex linear span of functions of $z$ from
$\{K^m(z,a):m=0,1,2,\dots\}$.  Similarly, define the
{\it Szeg\H o span\/} to be the complex linear span
$\mathcal S$ of all functions $h(z)$ of $z$ of the form $S^m(z,a)$
as $a$ ranges over $\Om$ and $m$ ranges over all non-negative integers,
and define ${\mathcal S}_U$ and ${\mathcal S}_a$ in the analogous
way.

We now list two theorems that will serve as the heart of the paper.
Note that if $f:\Om_1\to\Om_2$ is a biholomorphic mapping between
finitely connected domains in the plane, then it is well known that
there is a holomorphic branch of $\sqrt{f'}$ on $\Om_1$.

\begin{thm}
\label{thm3}
Suppose $\Om_1$ and $\Om_2$ are bounded domains in the plane bounded
by finitely many non-intersecting $C^\infty$-smooth Jordan curves.
Suppose further that $f:\Om_1\to\Om_2$ is a biholomorphic mapping.
Then $\Om_2$ is a boundary arc length quadrature domain
if and only if $f'$ is equal to the square of an element of the
Szeg\H o span associated to $\Om_1$, or equivalently, if and only if
$\sqrt{f'}$ is in the Szeg\H o span.
\end{thm}

We shall call a biholomorphic mapping $f$ as given in Theorem~\ref{thm3}
such that $\sqrt{f'}$ is in the Szeg\H o span a {\it Szeg\H o coordinate}.
This is a direct analogue of the definition of Bergman coordinate in \cite{B7}
where a biholomorphic map $f$ is called a Bergman coordinate if and only if
$f'$ is in the Bergman span.

The next theorem allows us to use Theorem~\ref{thm3} to approximate
domains by boundary arc length quadrature domains.

\begin{thm}
\label{thm4}
Suppose that $\Om$ is a bounded domain bounded by finitely many non-intersecting
$C^\infty$-smooth Jordan curves.  Then the Szeg\H o span $\mathcal S$ associated
to $\Om$ is dense in $A^\infty(\Om)$.  Furthermore, ${\mathcal S}_U$ is
dense in $A^\infty(\Om)$ for any open subset $U$ of $\Om$, and
${\mathcal S}_a$ is dense in $A^\infty(\Om)$ for any point $a\in\Om$.
\end{thm}

Theorem~\ref{thm4} was proved in \cite{B0}.

Theorems~\ref{thm3} and~\ref{thm4} and results of \cite{B7} and \cite{G2}
can be combined to obtain the following results.  The next theorem
shows that boundary arc length quadrature domains are dense in the realm of
bounded finitely connected smoothly bounded planar domains.

\begin{thm}
\label{thm5}
Suppose $\Om$ is a bounded domain bounded by finitely many non-intersecting
$C^\infty$-smooth Jordan curves.  There is a biholomorphic mapping $f$ which
is a Szeg\H o coordinate as close to the identity in $C^\infty(\Obar)$ as we
please.  Hence, $f(\Om)$ is a boundary arc length quadrature domain which is
$C^\infty$ close to $\Om$.
\end{thm}

The proof of Theorem~\ref{thm5} was sketched in \cite[p.~285-6]{B5}.  We shall
flesh out the argument here in \S\ref{sec3} to make this paper more
comprehensible and because we will need the fine details in order to
generalize the result to double quadrature domains.

\begin{thm}
\label{thm6}
Suppose $\Om$ is a bounded $n$-connected domain bounded by $n$ non-intersecting
$C^\infty$-smooth Jordan curves.  There is a biholomorphic mapping $f$ which is
both a Szeg\H o coordinate and a Bergman coordinate as close to the identity in
$C^\infty(\Obar)$ as we please.  Hence, $f(\Om)$ is both a boundary arc length
quadrature domain and an area quadrature domain, i.e., a double quadrature domain,
which is $C^\infty$ close to $\Om$.
\end{thm}

We will prove Theorem~\ref{thm6} in \S\ref{sec2} in case the domain is
simply connected and prove the general case in \S\ref{sec4}.

The Riemann mapping theorem says that every
simply connected domain can be mapped to a
special double quadrature domain, the unit disc.
Our results can be thought of as a generalization, or maybe
even an improvement.  Every finitely connected
domain in the plane can be mapped to a {\it nearby\/}
double quadrature domain.

Double quadrature domains have many remarkable
properties in common with the unit disc, some of which
we enumerate here.

$\bullet$ The boundary is given by finitely many non-intersecting real
algebraic curves which are $C^\infty$-smooth and real analytic.  (It is
shown in \cite{G} that the boundary is essentially given by the zero
set of an irreducible polynomial $P(z,w)$ on $\C^2$.  In fact, the
boundary is equal to $\{z: P(z,\bar z)=0\}$ minus perhaps finitely many
points, and $P(z,w)$ must have a a certain special form.)

$\bullet$ The complex polynomials belong to the Bergman span.

$\bullet$ The complex polynomials belong to the Szeg\H o span.  (In fact,
this property characterizes double quadrature domains.  See
Theorem~\ref{sequiv}.)

$\bullet$ The complex unit tangent vector function $T(z)$ is a rational
function of $z$ and $\bar z$ for $z$ in the boundary (which of course is
equivalent to being a rational function of $\text{Re }z$ and $\text{Im }z$).
The function $T(z)$ is also the restriction to the boundary of a meromorphic
function on the double without poles on the boundary.  This means that $T(z)$
can be extended to the domain as either a meromorphic function with no poles
on the boundary, or an anti-meromorphic function with no poles on the boundary.
(A more precise statement of this known fact can be found in
Theorem~\ref{thmX}.)

$\bullet$ The Schwarz function $S(z)$, which exists and extends to be analytic
on a neighborhood of the boundary by virtue of the real analyticity of the
boundary, also extends to be meromorphic on the domain.  It has no poles on the
boundary, and it is an algebraic function.  Because $S(z)=\bar z$ on the
boundary, both $z$ and $S(z)$ extend meromorphically to the double of the
domain.  It was noted in \cite{G} that all the meromorphic functions on the
double are generated by the extensions of $z$ and $S(z)$ to the double.

$\bullet$ Both the Bergman kernel $K(z,w)$ and the Szeg\H o kernel $S(z,w)$ are
rational functions of $z$, $S(z)$, $\bar w$, and $\overline{S(w)}$.  Hence,
they both extend to the double cross the double as functions which are
meromorphic in $z$ and anti-meromorphic in $w$.  If the domain is multiply
connected, the two kernels are rational functions of $f_1(z)$, $f_2(z)$,
$\overline{f_1(w)}$, and $\overline{f_2(w)}$, where $f_1$ and $f_2$ are
two Ahlfors maps associated to two generic points in the domain.  In the
simply connected case, the kernels are even more elementary (see below).

$\bullet$ Both the Bergman and Szeg\H o kernels are 
rational functions of $z$, $\bar z$, $w$, and $\bar w$ when $z$ and $w$
are restricted to the boundary.

$\bullet$ The Kerzman-Stein kernel $A(z,w)$ is a rational function of
$z$, $\bar z$, $w$, and $\bar w$ for $z$ and $w$ on the boundary.

If the domain is a simply connected double quadrature domain, then
even more can be said.  Let $f$ denote a Riemann map associated to
a point in the domain (which maps the domain one-to-one onto the
unit disc).

$\bullet$ $f$ extends to the double because $f(z)=1/\,\overline{f(z)}$
on the boundary.  Hence, it follows that $f(z)$ is a rational function of
$z$ and $S(z)$, and so it is a rational function of $z$ and $\bar z$ on
the boundary.  It is also algebraic.  Since the extension of $f$ to
the double generates the meromorphic functions on the double, it also
follows that $z$ and $S(z)$ are rational functions of $f(z)$.  (This
is another way to deduce that $f$ has a rational inverse, a well
know fact about simply connected quadrature domains with respect to area.)

$\bullet$
The Bergman kernel $K(z,w)$ and the Szeg\H o kernel $S(z,w)$
are rational functions of $f(z)$ and $\overline{f(w)}$.

$\bullet$
The Poisson kernel $p(z,w)$ is the real part of a function that is rational
in $z$ and $S(z)$, and $w$ and $\bar w$ (for $z$ in the domain and $w$ on
the boundary).  It is also the real part of a rational function of $f(z)$
and $w$ and $\bar w$.  It is even the real part of a rational function
of $f(z)$ and $f(w)$.

We shall explain these results in \S\ref{sec6}.  In \S\ref{sec7}, we
reveal a connection between the Schwarz function, quadrature
domains, and proper holomorphic mappings onto the unit disc.  We prove,
for example, that a domain is an area quadrature domain if and only if, on
the boundary, $z/\bar z$ is equal to the boundary values of a quotient
of proper holomorphic mappings to the unit disc (Theorem~\ref{aproper}).
A smooth domain is an arc length quadrature domain if and only if the
complex unit tangent vector function on the boundary is equal to the
boundary values of a quotient of proper holomorphic mappings to the
unit disc (Theorem~\ref{bproper}).

%

\section{Proof of Theorem~\ref{thm3}}
\label{sec5}
Suppose that $f:\Om_1\to\Om_2$ is a biholomorphic mapping between
bounded domains in the plane bounded by finitely many non-intersecting
$C^\infty$-smooth Jordan curves.  Let $F=f^{-1}$.  It is well known
that a holomorphic square root of $f'$ exists on $\Om_1$ and a square
root of $F'$ exists on $\Om_2$.  Since $f'$ and $F'$ extend
$C^\infty$-smoothly to the respective boundaries and are non-vanishing
there, the same is true of their square roots.  Notice that the change
of variables formula yields that
$$\int_{b\Om_1}|f'|(u\circ f)\ ds = \int_{b\Om_2}u\ ds$$
when $u$ is in the Hardy space $H^2(b\Om_2)$.  Let
$$\langle u,v\rangle_j=\int_{b\Om_j} u\,\bar v\ ds$$
denote the Hardy space inner product on $\Om_j$, $j=1,2$.
The last identity can be written
$$\langle \sqrt{f'}(u\circ f),\sqrt{f'}\rangle_1=
\langle u,1\rangle_2.$$
It is easy to verify that if $u\in L^2(b\Om_2)$ and $v\in L^2(b\Om_1)$,
then
$$\langle \sqrt{f'}(u\circ f),v\rangle_1=
\langle u,\sqrt{F'}(v\circ F)\rangle_2.$$

To prove the first half of Theorem~\ref{thm3}, assume that $\sqrt{f'}$
is in the Szeg\H o span associated to $\Om_1$.  If $h$ is an element
of the Hardy space $H^2(b\Om_2)$, then
$$\int_{b\Om_2} h\, ds = \int_{b\Om_1}|f'|(h\circ f)\ ds=
\langle \sqrt{f'}(h\circ f),\sqrt{f'}\rangle_1,$$
and if $\sqrt{f'}$ is in the Szeg\H o span, then this last inner product
yields a finite linear combination of values of $\sqrt{f'}(h\circ f)$ 
and its derivatives at finitely many points in $\Om_1$, which reduces
to a finite linear combination of values of $h$ and its derivatives at
finitely many points in $\Om_2$.  This shows that $\Om_2$ is a boundary
arc length quadrature domain.

To prove the converse, suppose that $\Om_2$ is a boundary arc length
quadrature domain.  Then, given any $g\in H^2(b\Om_1)$,
$$\langle g, \sqrt{f'}\rangle_1=\langle \sqrt{F'}(g\circ F), 1 \rangle_2,$$
and since $\Om_2$ is a boundary arc length quadrature domain, this last
inner product yields a finite linear combination of values of
$\sqrt{F'}(g\circ F)$ and its derivatives at finitely many points in
$\Om_2$, which reduces to a finite linear combination of values of $g$
and its derivatives at finitely many points in $\Om_1$.  Thus,
$\sqrt{f'}$ has the same effect as an element of the Szeg\H o span
when paired against a function $g$ in the Hardy space.  It follows that
$\sqrt{f'}$ must in fact be equal to that element in the Szeg\H o span.

\section{Density of quadrature domains in the simply connected case}
\label{sec2}
In this section, we give two proofs that boundary arc length quadrature
domains are dense in the realm of bounded $C^\infty$-smoothly bounded
{\it simply\/} connected domains, and that double quadrature domains are
dense as well.  The first proof introduces techniques that we will
generalize to the multiply connected setting.  The second proof is so
easy, natural, and short that one could get the wrong impression that
the multiply connected generalization should be easy and straightforward.

Suppose that $\Om$ is a bounded simply connected domain in the plane bounded
by a $C^\infty$-smooth Jordan curve.  We may now construct a boundary arc
length quadrature domain close to $\Om$ as follows.  By Theorem~\ref{thm4},
there exists a function $h$ in the Szeg\H o span that is close
to the constant function $1$ in $A^\infty(\Om)$.  (Of course, we take
$h$ close enough to $1$ so as to be non-vanishing.)  Let $H$ be a
holomorphic square root of $h$ on $\Om$, and let $f$ be a complex
antiderivative of $H$, where we choose the constant of integration so
that $f(z)$ is close to the identity in $A^\infty(\Om)$.  By
choosing $h$ close enough to $1$ in $A^\infty(\Om)$, we may guarantee
that $f$ is close enough to the identity that $f$ is one-to-one on
$\Om$ and $f(\Om)$ is as $C^\infty$-close as we desire to $\Om$.

Now according to Theorem~\ref{thm3}, $f(\Om)$ is a boundary arc
length quadrature domain that is $C^\infty$ close to $\Om$.

We can repeat the above argument using ${\mathcal S}_a$ for a point
$a$ in $\Om$ instead of the complete Szeg\H o span.  We obtain
a biholomorphic mapping $f$ to an arc length quadrature domain
close to $\Om$ satisfying an identity,
$$f'=\left( \sum_{m=0}^N c_m S^m(z,a) \right)^2.$$
We want to see that $f'$ is also in the Bergman span so that $f(\Om)$
is also an area quadrature domain via Theorem~1.3 of \cite{B7}.

Let $L(z,w)$ denote the Garabedian kernel associated to $\Om$ and
let $L^m(z,w)=(\dee/\dee w)^m L(z,w)$.  Let $\Lambda(z,w)$
denote the Schiffer kernel function, which is meromorphic in $z$ and $w$
with a double pole at $z=w$ and which is related to the Bergman
kernel $K(z,w)$ of $\Om$ via
$$K(w,z)\overline{T(z)}=-\Lambda(w,z)T(z)$$
when $w\in\Om$ and $z\in b\Om$.  Here $T(z)$ denotes the complex
unit tangent vector function pointing in the direction of the standard
orientation.  The basic properties of these kernel functions are
described in \cite{B1}.  We also define
$\Lambda^m(z,w)=(\dee/\dee w)^m \Lambda(z,w)$.  Let $b_j$ denote complex
coefficients that will be determined later.  The relationship
\begin{equation}
\label{SL}
S(w,z)=\frac{1}{i}L(z,w)T(z)
\end{equation}
for $w\in\Om$ and $z\in b\Om$ together with the relationship between $K(z,w)$
and $\Lambda(z,w)$ allow us to state that $H(z)T(z)=-\overline{G(z)T(z)}$ for
$z$ in $b\Om$, where
$$H(z)=\left[\left( \sum_{m=0}^N c_m S^m(z,a) \right)^2-
\left( \sum_{j=0}^M b_j K^j(z,a) \right)\right]$$
and
$$G(z)=\left[\left( \sum_{m=0}^N \bar c_m L^m(z,a) \right)^2-
\left( \sum_{j=0}^M \bar b_j \Lambda^j(z,a) \right)\right].$$
It can be deduced that $\left( \sum_{j=0}^N \bar c_j L^j(z,a) \right)^2$
is residue free since
$$\int_{b\Om}\left( \sum_{m=0}^N \bar c_m L^m(z,a) \right)^2\ dz=
\int_{b\Om} \left( \sum_{m=0}^N  c_m \overline{S^m(z,a)} \right)^2 \ d\bar z,$$
and this last integral is zero by Cauchy's Theorem.  Since
$\Lambda(z,a)$ has principal part $(-1/\pi)(z-a)^{-2}$ at $z=a$, it follows
that a positive integer $M$ and
constants $b_j$ can be chosen so that $G(z)$ has a removable singularity
at $a$.  Now both $H(z)$ and
$G(z)$ are holomorphic on $\Om$ and extend smoothly to the boundary in
such a way that $H(z)T(z)=-\overline{G(z)T(z)}$ on the boundary.  We now
claim that this can only be true if both $H$ and $G$ are identically zero.
Indeed, functions of the form
$\overline{G(z)T(z)}$ are orthogonal to holomorphic functions in the Hardy
space and functions of the form $H(z)T(z)$ are orthogonal to anti-holomorphic
functions.  Thus, $H(z)T(z)$ is orthogonal to holomorphic functions and
anti-holomorphic functions.  But every $C^\infty$-smooth function on $b\Om$
is the restriction to the boundary of a harmonic function in $C^\infty(\Obar)$,
and every such harmonic function can be written as $h+\overline{h}$ where
$h\in A^\infty(\Om)$.  Since smooth functions are dense in $L^2(b\Om)$,
it follows that $H\equiv0$ and $G\equiv0$.  We can now state that
$$\left( \sum_{m=0}^N c_m S^m(z,a) \right)^2=
\left( \sum_{j=0}^M b_j K^j(z,a) \right),$$
and this shows that $f'$ is in the Bergman span.  Hence $f$ is both a
Szeg\H o and a Bergman coordinate and $f(\Om)$ is a double quadrature
domain.

We remark here that, after finding the proof above, we discovered that
Avci \cite{Av} showed in his unpublished Stanford PhD thesis using other
methods that an arc length quadrature domain that satisfies a one point
quadrature identity would also have to be an area quadrature domain.
We include our proof here because it shows the close and explicit
connection between Szeg\H o and Bergman coordinates, and because we
will improve upon it in \S\ref{sec4} when we generalize our results
to smooth multiply connected domains.

We now turn to the simple proof of density of double quadrature
domains in the simply connected setting promised at the beginning
of this section.  The proof is based on a theorem of Aharonov and
Shapiro for area quadrature domains and a theorem of Shapiro and
Ullemar for arc length quadrature domains.  Let $\Om$ be a bounded
simply connected domain bounded by a $C^\infty$-smooth Jordan
curve and let $f$ denote a Riemann mapping of $\Om$ to the unit
disc.  Let $F=f^{-1}$.  Note that a branch of $\sqrt{F'}$ can be
chosen which is analytic on the unit disc and which extends to be
in $A^\infty(D_1(0))$.
Aharonov and Shapiro \cite{AS} proved that $\Om$
is an area quadrature domain if and only if $F$ is rational.
Shapiro and Ullemar \cite{SU} proved that $\Om$ is an arc
length quadrature domain if and only if $F'$ is the square
of a rational function.
To construct a double quadrature domain
that is $C^\infty$ close to $\Om$, we may approximate
$\sqrt{F'}$ in $A^\infty(D_1(0))$ by a polynomial $P(z)$.
Let $p(z)$ be (a polynomial) antiderivative of $P(z)^2$ where
the constant of integration is chosen so that $p(z)$ is
close to $F'$ in $A^\infty(D_1(0))$.  If the first polynomial
approximation $P(z)$ is sufficiently close to
$\sqrt{F'}$, then the image of the unit disc under $p(z)$
will be a double quadrature domain that is $C^\infty$ close
to $\Om$.

\section{Density of boundary arc length quadrature domains in the multiply
connected case}
\label{sec3}
In this section, we show how the argument in \cite[p.~77-78]{G2} using
half-order differentials and a special Runge Theorem on Riemann surfaces
can be modified along the lines of the previous section to prove
$C^\infty$ density of boundary arc length quadrature domains in the
category of smooth finitely connected domains.

Let $\Om$ be a bounded domain bounded by $n$ non-intersecting
$C^\infty$-smooth Jordan curves.  There is a function $h$ in the
Szeg\H o span that is as close to one in $C^\infty(\Obar)$ as we
desire.  We will need the rational functions $R_j(z)$, $j=1,\dots,n-1$
constructed in \cite[p.~78]{G2} that satisfy the following
special orthogonality conditions.  Let $\gamma_j$, $j=0,\dots,n-1$
denote the boundary curves of $\Om$ where $\gamma_0$ denotes the
outer boundary.  Choose a point $z_j$ inside $\gamma_j$,
$j=1,\dots,n-1$ (i.e., choose one point from each of the holes in
$\Om$).  The rational functions are given by
$$R_j(z)=\frac{1}{2\pi i(z-z_j)}+(z-z_j)Q_j(z),$$
where the $Q_j$ are any polynomials satisfying
$$Q_j(z_k)=-\frac{1}{2\pi i(z_k-z_j)^2}$$
for $k\ne j$ (making $R_j(z_k)=0$ for $k\ne j$).
The functions $R_j$ have the property that
$$\int_{\gamma_k} R_j\ dz = \delta_{kj}\quad\text{and}
\int_{\gamma_k} R_jR_m\ dz = 0$$
for $j,k,m=1,\dots n-1$.  Since each rational function $R_j$ is
in $A^\infty(\Om)$, we may approximate it in $A^\infty(\Om)$
by a function $r_j$ in the Szeg\H o span.  Notice that it
easy to find coefficients $A_j$ such that the periods of
$(h-\sum_{j=1}^{n-1}A_j R_j)^2$ vanish because the system is
$$\int_{\gamma_k} h^2\ dz + 2\sum_{j=1}^{n-1} A_j\int_{\gamma_k} hR_j\ dz=0$$
for $k=1,\dots,n-1$.  Since $h$ is close to one, the
coefficient matrix
$$\int_{\gamma_k} hR_j\ dz$$
is close to the identity and $\int_{\gamma_k} h^2\,dz$ is close to zero.
Thus, the $A_j$ are uniquely determined and are close to zero.
If we modify the problem of making the periods vanish by replacing
the $R_j$ by our functions $r_j$ in the Szeg\H o span, then the
system becomes
$$\int_{\gamma_k} h^2\ dz + 2\sum_{j=1}^{n-1} A_j\int_{\gamma_k} hr_j\ dz+
\sum_{j,m=1}^{n-1} A_jA_m\int_{\gamma_k} r_jr_m\ dz=0,$$
and since the coefficients in the linear part are close to the
identity matrix and the coefficients of the quadratic terms are
small, the implicit function theorem yields that the $A_j$ are
uniquely determined and are close to the small solution to the
system with $R_j$ in place of $r_j$.  Notice that by taking $h$
to be sufficiently close to one, the coefficients $A_j$ can be
made small enough that $h-\sum_{j-1}^{n-1}A_j r_j$ is as close
to one as desired.  Thus, this shows that there is a function
$F=h-\sum_{j-1}^{n-1}A_j r_j$ in
the Szeg\H o span which is close to one whose square $F^2$ has
vanishing periods. Now let $f$ be an antiderivative of $F^2$.
Since $F^2$ is close to one, we may choose the constant of integration
so that $f$ is close to the identity.  If $h$ is close enough to
one, then $f$ will be close enough to the identity to make $f$
biholomorphic on $\Om$ and Theorem~\ref{thm3} yields that $f$
maps $\Om$ to a boundary arc length quadrature domain that is
$C^\infty$ close to $\Om$.

\section{Density of double quadrature domains in the general case}
\label{sec4}
We continue with the same domain and setup as in the previous
section.  In particular, $\gamma_j$, $j=0,\dots,n-1$
denote the boundary curves of $\Om$ where $\gamma_0$ denotes the
outer boundary.  Let $p=n-1$, and let $\Ohat$ denote the double of
$\Om$.  We are now thinking of the curves $\gamma_j$ as being
curves on $\Ohat$.  Define $p$ more curves $\gamma_j$, $j=p+1,\dots,2p$,
on $\Ohat$ as follows.  Let $\Gamma_{j}$ be a simple smooth curve
that starts at a point on $\gamma_j$, enters $\Om$, and terminates
at a point on the outer boundary $\gamma_0$.  We may choose the
$\Gamma_{j}$ so that no two of them intersect.  In this case,
the domain $U$ given by $\Om-\cup_{j=1}^p\Gamma_{j}$ is a simply
connected domain.  Let $\widetilde{\Gamma_j}$ denote the curve
which is the image of the reverse of $\Gamma_{j}$ in the reflected
copy of $\Om$ in $\Ohat$, and let $\gamma_{p+j}$ be the closed
curve in $\Ohat$ obtained by following $\Gamma_{j}$ in $\Om$ and then
$\widetilde{\Gamma_j}$ in the reflected copy of $\Om$ back to the
starting point.  In this way, we obtain a standard homology
basis, $\gamma_j$, $j=1,\dots, 2p$, for $\Ohat$ where the first
$p$ curves are the standard homology basis for $\Om$ and the last
$p$ curves wrap around the ``handles'' of $\Ohat$ and complete
the homology basis for $\Om$ to a homology basis for $\Ohat$.

We will need the concept of a half-order differential on the
double as used in \cite{G2}.  In the present context, we will
describe a meromorphic half-order differential $h$ on $\Ohat$
as a pair of meromorphic functions $(h_1,h_2)$ where the two
functions are meromorphic on $\Om$, extend continuously to the
boundary (and are allowed to take the value $\infty$ at a boundary
point), and satisfy
\begin{equation}
\label{half}
h_1(z)T(z)=\overline{h_2(z)}
\end{equation}
on the boundary of $\Om$.  The symbol $h$ stands for the pair $(h_1,h_2)$.
These objects have the virtue that identity~(\ref{half}) allows the
meromorphic differential ${h_1}^2\,dz$ on $\Om$ to be extended to a meromorphic
differential on $\Ohat$ (via extension by $\overline{{h_2}^2}\,d\bar z$ on the
anti-holomorphic reflected side).  We let $h^2\,dz$ denote the extension of
$h_1^2\,dz$ to the double.  Another important feature is that if $h=(h_1,h_2)$
and $g=(g_1,g_2)$ represent two half-order differentials, then $h_1g_1\,dz$
extends meromorphically to $\Ohat$ as a meromorphic differential (via
extension by $\overline{h_2g_2}\,d\bar z$ on the reflected side).  We let
$hg\,dz$ denote the extension of this form to the double.  Hence, if $\gamma$
is a curve on $\Ohat$ and $h=(h_1,h_2)$ and $g=(g_1,g_2)$ are half-order
differentials, we may define the pairing
$$B(h,g)=\int_\gamma hg\ dz$$
as long as the extension $hg\,dz$ to $\Ohat$ does not have
poles along the curve.

Notice that identity~(\ref{half}) is the same
as the critical identity~(\ref{SL}) for the Szeg\H o and
Garabedian kernels.  Results in \cite{B2a} (Lemma~6.1)
show that all meromorphic half-order differentials
without poles on the boundary arise as linear combinations of
the Szeg\H o and Garabedian kernels and their derivatives
in the second variable.  When there are no poles inside
the domain, a meromorphic half-order differential must be
a linear combination of just the Szeg\H o kernel and its derivatives
in the second variable on the domain side.  This is the key
connection between the Szeg\H o kernel and the results from
\cite{G2} that allow us to use Szeg\H o coordinates to
prove the density of double quadrature domains here.

The idea of the proof is to construct a half-order differential
$h=(h_1,h_2)$ with the property that $h_1$ is in the Szeg\H o
span and is close to one in $A^\infty(\Om)$, and $h^2 dz$ is a
meromorphic one-form on $\Ohat$ with a single pole at a
point in the reflected side of $\Om$ in $\Ohat$ with
the property that $\int_{\gamma_j} h^2\,dz=0$ for $j=1,\dots,2p$.
Then $h^2\,dz=df$ where $f$ is a meromorphic function on $\Ohat$.
If $h_1$ is sufficiently close to one, then the condition
$f'={h_1}^2$ implies that $f$ is a biholomorphic mapping of
$\Om$ onto a domain which is $C^\infty$ close to $\Om$.
Now Theorem~\ref{thm3} yields that $f(\Om)$ is a boundary arc
length quadrature domain.  Finally, by \cite{G}, the fact
that $f$ extends meromorphically to the double yields that
$f(\Om)$ is an area quadrature domain.  In this way, we will
obtain a double quadrature domain that is $C^\infty$ close
to $\Om$.

We shall need the following special Mergelyan approximation
lemma.  Suppose $a_0\in U$.  Recall that ${\mathcal S}_{a_0}$
denotes the complex linear span of functions of $z$ from
$\{S^m(z,a_0):m=0,1,2,\dots\}$.

\begin{lem}
\label{mergelyan}
Given a holomorphic function $G$ in $A^\infty(\Om)$ and
continuous functions $\varphi_j$ on $\Gamma_j$ that satisfy the
compatibility condition
$$G(w)T(w)=\overline{\varphi_j(w)}$$
at the endpoints $w$ of each curve $\Gamma_j$ (which fall on
$b\Om$), it is possible to find a function
$$\sigma(z)=\sum_{m=0}^Nc_mS^m(z,a_0)$$
in ${\mathcal S}_{a_0}$ such that
$\sigma$ is as close as desired to $G$ in $A^\infty(\Om)$ and
so that 
$$\lambda(z)=-i\sum_{m=0}^N\overline{c_m}L^m(z,a_0)$$
is as close to $\varphi_j$ on each $\Gamma_j$ in the uniform
topology as desired.  (Note that identity (\ref{SL}) shows
that $(\sigma,\lambda)$ is a half-order differential.)
\end{lem}

We will now prove the density of double quadrature domains,
assuming the truth of the approximation lemma.  Afterwards,
we will prove Lemma~\ref{mergelyan}.

Let $K$ denote the compact set
$\Obar\cup\left(\cup_{j=1}^p\gamma_{j+p}\right)$
in $\Ohat$.  We may consider non-meromorphic half-order
differentials defined on $K$ where $g=(g_1,g_2)$ is such
that $g_1$ is in $C^\infty(\Obar)$ and the $g_2$ is a continuous
function on the union of the curves $\Gamma_j$ and satisfies the
compatibility condition
$$g_1(w)T(w)=\overline{g_2(w)}$$
at the endpoints $w$ of each curve $\Gamma_j$.
For two such half-order differentials $g$ and $h$, we may define
the pairing
$$B_k(g,h)=\int_{\gamma_k} gh\ dz$$
for $k=1,\dots,2p$ as described above in the meromorphic case.

We now define $2p+1$ such half-order differentials $g^m=(g_1^m,g_2^m)$,
$m=0,\dots,2p$, on $K$ satisfying

\begin{enumerate}[(i)]
\item $g_1^0=1$ on $\Obar$,
\label{i}
\item $B_k(g^0,g^0)=0$ for all $k=1,\dots,2p$, and
\label{ii}
\item $B_k(g^0,g^m)=\delta_{km}$ for $k,m=1,\dots,2p$,
\label{iii}
\end{enumerate}

\noindent
We start by setting $g_1^0\equiv1$ on $\Obar$
(where $g_1^0$ is the first function associated to the half-order
differential $g^0$) as condition (\ref{i}) dictates.  Notice that
(\ref{ii}) is therefore automatically satisfied for $k=1,\dots,p$.
It is now an easy matter to define $g_2^0$ on the curves $\Gamma_j$
so that the compatibility conditions hold and condition (\ref{ii}) holds
for $k=p+1,\dots,2p$ as well.
Next, set $g_1^m=R_m$ on $\Obar$, where $R_m$ are the rational
functions of the previous section for $m=1,\dots,p$.  Note that
condition (\ref{iii}) now holds for $k,m=1,\dots,p$.  It is another
simple matter to define $g_2^m$ on the curves $\Gamma_j$ so that
(\ref{iii}) also holds for $k=p+1,\dots,2p$ and $m=1,\dots,p$.
Finally, $g^m$ are chosen for $m=p+1,\dots,2p$ so that $g_1^m\equiv0$
on $\Obar$ and $g_2^m$ are defined on the curves $\Gamma_j$ so that
(\ref{iii}) holds for $k=p+1,\dots,2p$.  Note that (\ref{iii})
holds automatically for $k=1,\dots,p$.

Lemma~\ref{mergelyan} implies that there are half-order differentials
$h^m=(h_1^m,h_2^m)$ which are close to the $g^m$ in the sense that
$h_1^m$ is close to $g_1^m$ in $A^\infty(\Om)$ for each $m$ and
$h_2^m$ is uniformly close to $g_2^m$ on the curves $\Gamma_j$.
The half-order differential that we seek will be given by
$$h=h^0+\sum_{m=1}^{2p}t_mh^m,$$
where the $t_m$ are complex numbers close to zero.
(Note that we may take linear combinations of half-order
differentials in the obvious way by summing the corresponding
functions in the pairs.)  We will choose the $t_m$ to make
$B_k(h,h)=0$ for $k=1,\dots,2p$ so that $h^2\, dz$ is period free,
and is therefore an exact meromorphic differential $df$.  (Note that
the lone pole of $h^2\,dz$ would necessarily have residue equal to zero.)
Let $g=g^0+\sum_{m=1}^{2p}t_mg^m$.
Define complex coefficients via
\begin{eqnarray*}
a_{kij} & = & B_k(h^i+g^i,h^j-g^j) \\
b_{kj} & = & 2\left(B_k(g^j,h^0-g^0)+B_k(h^0,h^j-g^j)\right) \\
c_k & = &B_k(h^0+g^0,h^0-g^0) \\
A_{kij} & = & B_k(g^i,g^j)
\end{eqnarray*}
Note that $A_{kij}$ is fixed, while the others depend on the
$h^j$ and can be made arbitrarily small by improving the
approximation sufficiently.
Expand $B_k(h,h)=B_k(h-g+g,h-g+g)$ to obtain
$$B_k(h,h)=2t_k+c_k+\sum_{i,j}A_{kij}t_it_j
+\sum_{i,j} a_{kij}t_it_j
+\sum_{j} b_{kj}t_j.$$
We obtain a system by setting these $2p$ periods equal
to zero.  Note that this system has the unique zero solution if
the coefficients that depend on the $h^j$ are taken to be
zero.  It is now a standard application of the implicit function
theorem to deduce that it is possible to choose $t_j$
close to zero to make all these periods vanish.
This completes the proof, assuming Lemma~\ref{mergelyan}.
We now turn to proving the lemma.

Note that Theorem~\ref{thm4} allows us to simplify the assumptions
in Lemma~\ref{mergelyan} by subtracting off an element of the
Szeg\H o span to be able to assume that $G$ is {\it zero} on $\Obar$.
Given the point $a_0$ in the simply connected domain
$U=\Om-\cup_{j=1}^p\Gamma_{j}$, our
task then is to find a function $\lambda(z)$ which is a linear
combination of functions of the form $L^m(z,a_0)$ which is
uniformly close to the functions $\varphi_j$ on each $\Gamma_j$
and which is close in $C^\infty(b\Om)$ to zero.  This turns
out to be remarkably easy.
The first step in the argument is to use the classical Mergelyan
Theorem (see Rudin \cite{R}, Chapter~20, Exercise~1) and Runge's
Theorem to approximate our functions as follows.  For $\epsilon>0$, let
$V_\epsilon$ denote the points in $\mathbb C$ that are less than or equal
to a distance $\epsilon$ from $b\Om$, i.e., the closure of a small collared
neighborhood of the boundary in $\mathbb C$.
Let $K_\epsilon =V_\epsilon\cup\left(\cup_{j=1}^p\Gamma_{j}\right)$.
We first use Mergelyan's Theorem to obtain a rational
function $R(z)$ with poles in ${\mathbb C}-K_\epsilon$ that is
unformly close to zero in $V_\epsilon$, and because the functions
$\varphi_j$ vanish on $b\Om$, and because we may shrink $\epsilon$,
we may also arrange for $R(z)$ to be
uniformly close to the $\varphi_j$ on each $\Gamma_j$.  Cauchy's
Estimates reveal that uniform convergence on the interior of
$V_\epsilon$ implies $C^\infty$ convergence on $b\Om$.  Hence, we
may obtain a rational function that is close to zero in $C^\infty(b\Om)$
and uniformly close to $\varphi_j$ on each $\Gamma_j$.  Next,
we may slide the poles of $R(z)$ that fall in $U$ to the point
$a_0$ as is done in many standard proofs of Runge's Theorem (see
Stein \cite[p.~63]{St}) to obtain a rational function $R(z)$
with poles at $a_0$ and finitely many other points outside of $\Obar$
which is as close in $C^\infty(b\Om)$ to zero as desired and as
close in the uniform topology to $\varphi_j$ on each $\Gamma_j$
as desired.

The Garabedian kernel $L(z,w)$ is such that
$$L(z,w)=\frac{1}{2\pi(z-w)}+\ell(z,w),$$
where $\ell(z,w)$ is holomorphic in $z$ and $w$ and is in
$C^\infty(\Obar\times\Obar)$ (see \cite[p.~102]{B1}).  We
will therefore be able to use $L(z,w)$ much like a Cauchy
kernel in what follows.
The Residue Theorem yields that, for $z\in\Om-\{a_0\}$,
$$\frac{1}{i}\int_{w\in b\Om}L(z,w)R(w)dw$$
is equal to $R(z)$ plus $2\pi$ times the residue of
$L(z,w)R(w)$ at $a_0$, which is a fixed finite linear combination
of functions of the form $L^m(z,a_0)$.  If we show that
the integral tends to zero in $z$ in $A^\infty(\Om)$, then
we will have completed the proof.  This turns out to be an
easy consequence of the fact that the Szeg\H o projection
is a continuous linear operator in $C^\infty(\Obar)$ (see
\cite[p.~13]{B1}).  Indeed, identity (\ref{SL}) and the
fact that $L(z,w)=-L(w,z)$ yields
that the integral is equal to
$$-\int_{w\in b\Om}S(z,w)R(w)ds,$$
which is the Szeg\H o projection evaluated at $z$ of the
function which is the restriction of $R(w)$ to the boundary.
Since $R$ is close to zero in $C^\infty(b\Om)$, it follows
that the function given by the integral is close to zero
in $A^\infty(b\Om)$.  The proof of the lemma is finished.

\section{Properties of double quadrature domains}
\label{sec6}
Suppose $\Om$ is a double quadrature domain in the plane.
Because $\Om$ is an area quadrature domain, it has a boundary
given by an algebraic curve, and because it is a boundary arc
length quadrature domain, that curve consists of $C^\infty$ smooth
real analytic curves (see \cite{AS, G, G2}).  Because $\Om$ is
an area quadrature domain, the complex polynomials belong to the
Bergman span (see \cite{B7}).  Also, the Schwarz function $S(z)$ 
exists and is meromorphic on $\Om$, extends analytically past
the boundary, and satisfies $S(z)=\bar z$ on the boundary (see
\cite{AS,G}).  This identity shows that the two functions $z$ and
$S(z)$ extend meromorphically to the double of $\Om$.  It was proved
in \cite{G} that the extensions form a primitive pair for the field
of meromorphic functions on the double, i.e., all the meromorphic
functions on the double are given as rational combinations of
the two.  It is shown in \cite{B5} (see also \cite{B4}) that
the Bergman kernel $K(z,w)$ associated to an area quadrature
domain is a rational function of $z$, $S(z)$, and
$\bar w$, and $\overline{S(w)}$.  Consequently, $K(z,w)$ is
a rational function of $z$, $\bar z$, $w$, and $\bar w$ when
$z$ and $w$ are restricted to the boundary.  It was also
shown in \cite{B5,B4} that the Bergman kernel is a rational
combination of two Ahlfors maps associated to two generic
points in the domain.  Aharonov and Shapiro \cite{AS} showed
that the Ahlfors maps are algebraic.

Because $\Om$ is an arc length quadrature domain, there is a meromorphic
function $H(z)$ on $\Om$ which extends analytically past the boundary
and such that $\overline{H(z)}$ is equal to the complex unit tangent
vector function $T(z)$ on the boundary (see \cite{G2,Av}).  The
identity relating the Szeg\H o and Garabedian kernels can be used
with this fact to see that both the Szeg\H o and the Garabedian
kernels extend to the double in the first variable when the second
variable is held fixed.  Indeed, the identity
$$\overline{S(z,a)}=\frac{1}{i}L(z,a)T(z)\qquad\text{for }z\in b\Om$$
shows that $L(z,a)=i\overline{S(z,a)}/\,\overline{H(z)}$ on the boundary,
and this shows that $L(z,a)$ extends meromorphically to the double in $z$.
Similarly, $S(z,a)$ extends.  These identities also reveal that
$T(z)$ is equal to a function $G(z)$ on the boundary where $G(z)$ is
meromorphic on $\Om$ and extends analytically past the boundary.
Hence, $T(z)$ is the restriction to the boundary of a meromorphic
function on the double of $\Om$ that has no singularities on the
boundary of $\Om$.  These arguments are reversible and we may
state the following theorem.

\begin{thm}
\label{thmX}
A bounded domain with a piecewise $C^1$ smooth boundary is a
boundary arc length quadrature domain if and only if the Szeg\H o
kernel $S(z,w)$ associated to the domain extends meromorphically
to the double as a function of $z$ for each fixed $w$ in the domain.
This is the case if and only if the complex unit tangent vector
function is the restriction to the boundary of a meromorphic
function that has no singularities on the boundary of $\Om$.
\end{thm}

The same theorem holds with the Garabedian kernel in place
of the Szeg\H o kernel.  (The part of the theorem about $T$
was proved in \cite{G2}.  The equivalence of the extendibility
of the Szeg\H o kernel is a rather direct consequence of the
extendibility of $T$ and might very well have been noted earlier.)

If $\Om$ is a double quadrature domain, then we have shown
that $T(z)$ is the restriction to the boundary of a meromorphic
function on the double, and we also know that such functions
are generated by $z$ and $S(z)$.  Hence, it follows that $T(z)$
is the restriction of a rational combination of $z$ and $S(z)$ to the
boundary, i.e., $T(z)$ is a rational function of $z$ and $\bar z$.

When Theorem~\ref{thmX} is combined with results from \cite{B2},
it follows that the Szeg\H o kernel $S(z,w)$ associated to
a double quadrature domain is a rational combination of
$z$, $S(z)$, and $\bar w$, and $\overline{S(w)}$.  Consequently,
$S(z,w)$ is a rational function of $z$, $\bar z$, $w$, and $\bar w$
when $z$ and $w$ are restricted to the boundary.

The Kerzman-Stein kernel is given by
$$A(z,w)=\frac{1}{2\pi i}
\left(\frac{T(w)}{w-z}-\frac{\overline{T(z)}}{\bar w-\bar z}\right)$$
for $z$ and $w$ in the boundary of $\Om$, and it follows that
the Kerzman-Stein kernel associated to a double quadrature domain
is a rational function of $z$, $\bar z$, $w$, and $\bar w$.

We shall now show that the complex polynomials are in the Szeg\H o
span associated to a double quadrature domain.  Let $H(z)$ denote
the meromorphic function that is equal to $\overline{T(z)}$ on the
boundary.  Suppose $h\in A^\infty(\Om)$.  Since $\bar z=S(z)$ on the
boundary, we may write
$$\int_{b\Om} h\, {\bar z}^n\ ds = \int_{b\Om} h(z)S(z)^n\, \overline{T(z)}\,dz =
\int_{b\Om} h(z)S(z)^n H(z)\,dz,$$
and this last integral is equal to a finite linear combination of
values of $h$ and its derivatives at finitely many points in $\Om$ by
the Residue Theorem.  Hence, $z^n$ has the same effect when paired with
$h$ that a certain element in the Szeg\H o span would have.  Since
$A^\infty(\Om)$ is dense in the Hardy space, it follows that $z^n$
is equal to the element in the Szeg\H o span.

The converse of this last result is true, namely, that if the complex
polynomials are in the Szeg\H o span, then the domain is a double
quadrature domain.  To see this, note that $z=z/1$.  If both $z$ and
$1$ are in the Szeg\H o span, then $z$ is a quotient of elements in the
Szeg\H o span, and as such, it extends to the double as a meromorphic
function (since the relationship between the Szeg\H o and the Garabedian
kernels reveals that quotients of functions in the Szeg\H o span would
be equal to quotients of conjugates of functions in the ``Garabedian
span'' on the boundary).  The condition that $z$ extends to the double
is the classic condition that is equivalent to the domain being an area
quadrature domain.  It is easy to see that a smooth domain is a boundary
arc length quadrature domain if and only if the constant function $1$ is
in the Szeg\H o span.  Hence, we have proved the following theorem.

\begin{thm}
\label{sequiv}
A bounded domain with a piecewise $C^1$ smooth boundary is a double
quadrature domain if and only if the complex polynomials belong to
the Szeg\H o span associated to the domain.
\end{thm}

Theorem~\ref{sequiv} is an appealing analogue of Theorem~1.2 from
\cite{B7} which states that a bounded domain is an area quadrature
domain if and only if the complex polynomials belong to the Bergman
span.

The Poisson kernel of a smooth simply connected domain is given by
$$p(z,w)=
2\text{Re } \left(\frac{S(z,w)S(w,a)}{S(a,a)}\right)-\frac{|S(w,a)|^2}{S(a,a)},$$
where $z\in\Om$, $w\in b\Om$, and $a$ is a fixed point in $\Om$ (see
\cite[p.~1367]{B}).  Hence, when $\Om$ is a simply connected double quadrature
domain, $p(z,w)$ is the real part of a function that is rational in $z$ and
$S(z)$, and $w$ and $\bar w$ (just like in the unit disc, where $S(z)=1/z$).

The rest of the properties mentioned in \S1 involving the Riemann map $f$
of a simply connected domain
follow from the remarks above and the fact that $f$ extends meromorphically
to the double and the extension generates the field of meromorphic functions
on the double, i.e., every meromorphic function on the double is a rational
function of the extension of the Riemann map.  Hence, for example, on a
quadrature domain with respect to area, the functions $z$ and $S(z)$ extend
to the double, and so they are rational functions of $f(z)$.  Since,
$f(w)=1/\,\overline{f(w)}$ on the boundary, the restriction to the boundary
of a rational function of $f(w)$ and $\overline{f(w)}$ is equal to a rational
function of $f(w)$ alone.  Hence, it follows from the remarks above about the
Poisson kernel is the real part of a rational function of $f(z)$ and $f(w)$.

\section{The Schwarz function and proper holomorphic mappings to the disc}
\label{sec7}
We conclude this paper by noting a relationship between
proper holomorphic mappings of a domain to the unit disc and the
existence of the Schwarz function, which is the hallmark of a
bounded quadrature domain with respect to area.  We go on to reveal
a further relationship between proper holomorphic mappings and
the extension of the unit tangent function to the double, which
characterizes a bounded quadrature domain with respect to area.

Assume that $\Om$ is a (bounded) area quadrature domain.
For the moment, assume that zero does not belong to $\Om$.
As remarked earlier, the Schwarz function $S(z)$ exists
and is meromorphic in $\Om$.  Let
$\{a_j\}_{j=1}^N$ denote the poles of $S(z)$ in $\Om$, and suppose
$a_j$ is a pole of $S(z)$ of order $n_j$.  Let
$f_j$ denote the Ahlfors map associated to $a_j$, which is a proper
holomorphic mapping of $\Om$ onto the unit disc such that $a_j$ is
a simple zero of $f_j$.  (The Ahlfors map associated to a point $a$
in $\Om$ is the unique solution to the extremal problem to maximize
$h'(a)$ under the conditions that $h$ maps $\Om$ into the unit disc
and $h'(a)$ is real.  It maps a bounded $n$-connected domain onto the
unit disc as an $n$-to-one branched covering map.)  We note here
that, since the boundary of $\Om$ is piecewise smooth, proper
holomorphic mappings to the unit disc extend continuously to the
boundary.  Conversely, a holomorphic
mapping to the unit disc that is continuous up to the boundary is
proper if and only if it is non-constant and it maps the boundary
to the unit circle.  Notice that,
since $|f_j(z)|=1$ and $\bar z=S(z)$ on the boundary, the function
$$\Phi(z)=\left(\prod_{j=1}^N f_j(z)^{n_j}\right)\frac{S(z)}{z}$$
has removable singularities at each $a_j$, is continuous up to the
boundary, and has unit modulus on the boundary.  If $\Om$ is multiply
connected, then $\Phi$ must have a zero in $\Om$ since the Schwarz
function must have a zero in the domain (see \cite{Av}).  Hence
$\Phi$ is a proper holomorphic mapping of $\Om$
to the unit disc.  It follows that $S(z)$ is $z$ times the quotient of
two proper holomorphic mappings of $\Om$ to the unit disc.  If the
point zero were in the domain, we could add the Ahlfors map $f_0$
to the product and obtain the same result.  If $\Om$ is simply
connected, then we can define
$$\Psi(z)=\left(\prod_{j=1}^N f_j(z)^{n_j}\right)
\left(\prod_{k=1}^M F_k(z)^{-m_k}\right)\frac{S(z)}{z},$$
where $f_j$ is the Riemann map associated to a point $a_j$
where $S(z)$ has a pole of order $n_j$, and $F_k$ is the Riemann
map associated to a point $b_k$ where $S(z)$ has a zero of
order $m_k$.  In this case, $\Psi$ has constant modulus one
on the boundary and no zeros in $\Om$.  Hence it is a unimodular
constant and we deduce that
$S(z)$ is equal to a unimodular constant times $z$ times a
quotient of products of Riemann maps.  Since proper maps are
all given by unimodular constants times products of Riemann
maps in this setting, the same result follows, i.e., that
$S(z)$ is $z$ times the quotient of two proper holomorphic
mappings of $\Om$ to the unit disc.

If the domain $\Om$ is bounded by finitely many non-intersecting
Jordan curves, a simple converse follows.  Indeed, since a proper
holomorphic mapping $f$ to the unit disc would be continuous up
to the boundary and have unit modulus on the boundary in this case,
the identity $f(z)=1/\,\overline{f(z)}$, which holds on the boundary,
reveals that $f$ extends to the double.  If $\bar z$ is equal to $z$
times the quotient of two proper holomorphic mappings to the unit disc
on the boundary, then $z$ is seen to extend to the double, and it
follows that the domain is an area quadrature domain.  The theorem with
its converse reads as follows.

\begin{thm}
\label{aproper}
A bounded domain in the plane bounded by finitely many
non-intersecting Jordan curves is an area quadrature domain
if and only if, on the boundary, the function $\bar z/z$ is
equal to the boundary values of a quotient of proper holomorphic
mappings of the domain onto the unit disc.
\end{thm}

Note that if $\bar z/z=f_1/f_2$ on the boundary, then
the Schwarz function is given by $S(z)=zf_1(z)/f_2(z)$.

If a domain has piecewise $C^1$ smooth boundary, it is a
boundary arc length quadrature domain if and only if the
unimodular function $T(z)$ extends to the double as a
meromorphic function.  Similar reasoning to the argument
above with $T(z)$ in place of $S(z)/z$ yields the following
result.

\begin{thm}
\label{bproper}
A finitely connected bounded domain with piecewise $C^1$
smooth boundary is a boundary arc length quadrature
domain if and only if the complex unit tangent vector function
$T(z)$ is equal to a quotient of proper holomorphic mappings
of the domain onto the unit disc.
\end{thm}

Since an area quadrature domain has piecewise $C^1$
smooth boundary, we may also state this last theorem.

\begin{thm}
\label{dproper}
A bounded domain in the plane bounded by finitely many
non-intersecting Jordan curves is a double quadrature domain
if and only if, on the boundary, the function $\bar z/z$ is
equal to the boundary values of a quotient of proper holomorphic
mappings of the domain onto the unit disc, and the complex
unit tangent vector function $T(z)$ is also
equal to the boundary values of a quotient of proper holomorphic
mappings of the domain onto the unit disc.
\end{thm}

Finally, we remark that the semi-group of all proper holomorphic
mappings of a finitely connected domain to the unit disc has been
described completely in \cite{BK}.

\end{document}